\documentclass[11pt]{article}
\usepackage{amsfonts,amssymb,amsmath,bm,amsthm}
\usepackage{mathrsfs}
\newcommand{\Badr}{\mathbf{Bad}(\rr)}
\newtheorem{theorem}{Theorem}
\newtheorem{corollary}{Corollary}

\theoremstyle{remark}
\newtheorem{remark}{Remark}
\theoremstyle{definition}

\theoremstyle{definition}
\newtheorem{problem}{Problem}

\def\R{\mathbb{R}}

\def\Q{\mathbb{Q}}
\def\Z{\mathbb{Z}}
\def\N{\mathbb{N}}

\def\cA{\mathcal{A}}
\def\cU{\mathcal{U}}
\def\cK{\mathcal{K}}
\def\cW{\mathcal{W}}

\def\cH{\mathcal{H}}
\def\cB{\mathcal{B}}
\def\cM{\mathcal{M}}

\def\cC{\mathcal{C}}

\def\cF{\mathcal{F}}

\newcommand{\SL}{{\rm SL}}

\newcommand{\ve}{\varepsilon}
\newcommand{\vv}[1]{{\mathbf{#1}}}

\newcommand{\codim}{\operatorname{codim}}

\newcommand{\GL}{\mathrm{GL}}

\newcommand{\V}{\R^{n+1}}

\newcommand{\dt}{\cdot}
\newcommand{\we}{\wedge}
\def\Bad{\mathbf{Bad}}
\newcommand{\rr}{\mathbf{r}}

\newcommand{\da}{Diophantine approximation}

\usepackage{a4wide}

\usepackage{xcolor}

\numberwithin{equation}{section}

\begin{document}

\large

\title{\bf Quantitative non-divergence and\\ Diophantine approximation on manifolds}

\author{V.~Beresnevich\footnote{ORCID: 0000-0002-1811-9697 (also known as V. Berasnevich).} \and D.~Kleinbock\footnote{Supported in part by NSF grant DMS-1600814.}}

\date{{\small \em Dedicated to G. A. Margulis}}%

\maketitle

\vspace*{-3ex}

\begin{abstract}
The goal of this survey is to discuss the Quantitative non-Divergence estimate on the space of lattices and present a selection of its applications. The topics covered include extremal manifolds, Khintchine-Groshev type theorems, rational points lying close to manifolds and badly approximable points on manifolds. The main emphasis %of the chapter
 is on the role of the Quantitative non-Divergence estimate in the aforementioned topics within the theory of Diophantine approximation, and therefore this paper should not be regarded as a comprehensive overview of the area.
\end{abstract}

\section{Quantitative non-Divergence estimate and its origins}

\subsection{Background}\label{Back}

The main purpose of this survey is to discuss a particular strand of fruitful interactions between Diophantine approximation and the methods of Homogeneous Dynamics. The focus will be on the technique/estimate developed in \cite{MR1652916} by Margulis and the second named author %of this chapter
which is commonly known by the name of {\em Quantitative non-Divergence} ({\em abbr.} QnD). Before considering any quantitative aspects of the theory, it will be useful to explain the meaning of `{\em non-divergence}' of sequences and maps in the space
$$
X_k:=\SL_k(\R)/\SL_k(\Z)
$$
of real unimodular lattices. As is well known, {the quotient topology induced from $\SL_k(\R)$ makes this space non-compact}. Naturally, a non-divergent sequence in $X_k$ is then defined by requiring that it keeps returning into %a sufficiently large
{some compact} set. To give this narrative description more rigour it is convenient to use Mahler's compactness theorem and the function
$$
\delta:X_k\to\R_+
$$
which assigns the length of the shortest non-zero vector to a given lattice. Thus,
$$
\delta(\Lambda):=\inf\big\{\|\vv v\|:\vv v\in\Lambda\smallsetminus\{\vv0\}\big\}\quad\text{for every}\quad \Lambda\in X_k\,.
$$
Mahler's Compactness Theorem \cite{MR0017753} states that {\em a subset $S$ of $X_k$ is {relatively compact}
%bounded
if and only if there exists $\ve>0$ such that $\delta(\Lambda)\ge\ve$ for all $\Lambda\in S$.}
Thus, a sequence of lattices is non-divergent if and only if for a suitably chosen $\ve>0$ the sequence contains infinitely many elements in the (compact) set
\begin{equation}\label{K_e}
\cK_\ve:=\big\{\Lambda\in X_k:\delta(\Lambda)\ge\ve\big\}.
\end{equation}
The choice of the norm $\|\cdot\|$ does not affect Mahler's theorem. For simplicity we shall stick to the supremum norm: $\|\vv v\|=\max_{1\le i\le k}|v_i|$ for $\vv v=(v_1,\dots,v_k)$.

Similarly, given a continuous map
$$
\phi:[0,+\infty)\to X_k\,,
$$
we will say that $\phi(x)$ is {\em non-divergent} (as $x\to+\infty$) if there exists $\ve>0$ such that
$
\phi(x)\in\cK_\ve
$
for arbitrarily large $x$.

The development of the QnD estimate in \cite{MR1652916} was preceded by several important non-quantitative results instigated by Margulis \cite{MR0291352} regarding the orbits of one-parameter unipotent flows. The main result of \cite{MR0291352} verifies that if $\{u_x\}_{x\in\R}$ is a one-parameter subgroup of $\SL_k(\R)$ {consisting} of unipotent matrices, then $\phi(x)=u_x\Lambda$ is non-divergent for any $\Lambda\in X_k$. Several years later Dani \cite{MR530631} strengthened Margulis' result by showing that such orbits return into a suitably chosen compact set with positive frequency. To be more precise, Dani proved that there are $0<\ve,\eta<1$ such that for any interval $[0,t]\subset[0,+\infty)$ one has that
\begin{equation}\label{v1}
\lambda\big\{x\in[0,t]:\phi(x)\not\in\cK_\ve\big\}<\eta\,t\,,
\end{equation}
where $\lambda$ stands for Lebesgue measure on $\R$.
Subsequently Dani \cite{MR857195} improved his result by showing that under a {mild additional constraint} on $\{u_x\}_{x\in\R}$ the parameter $\eta>0$ {in} \eqref{v1} can be made arbitrarily small, in which case, of course, $\ve$ has to be chosen appropriately small. Later Shah \cite{MR1291701} generalised Dani's result to polynomial maps $\phi$ that are not necessarily orbits of some subgroups of $\SL_k(\R)$. It has to be noted that {the non-divergence theorem of Margulis was used as an ingredient in his proof of arithmeticity of non-uniform lattices in semisimple Lie groups of higher rank \cite{MR0422499}, and that subsequent qualitative non-divergence estimates, in particular, Dani's result in \cite{MR857195}}, were an important part of various significant developments of the time such as Ratner's celebrated theorems \cite{MR1262705}. The essence of the Quantitative non-Divergence estimate obtained in \cite{MR1652916} is basically an explicit dependence of $\eta$ on $\ve$ {in} \eqref{v1}. More to the point, it is applicable to a very general class of maps $\phi$ of several variables that do not have to be polynomial, let alone the orbits of unipotent subgroups. In the next subsection we give the precise formulation of the QnD estimate.

\subsection{The Quantitative non-Divergence estimate}

First we recall some notation and definitions. Given a ball $B={B}(x_0,r)\subset\R^d$ centered at $x_0$ of radius $r$ and $c>0$, {by}~$cB$ we will denote the ball $B(x_0,cr)$. Throughout, $\lambda_d$ will denote Lebesgue measure on $\R^d$. Given an open subset $U\subset \R^d$ and real numbers $C,\alpha>0$, a function $f:U\to\R$ is called {\em $(C,\alpha)$-good on $U$}\/ if for any ball $B\subset U$ for any $\ve>0$
$$
\lambda_d\big(\big\{x\in B:|f(x)|<\ve\big\}\big)~\le~ C\left(\frac{\ve}{\sup_{x\in B}|f(x)|}\right)^\alpha\lambda_d(B)\,.
$$
Finally, given $\vv v_1,\dots,\vv v_r\in\R^k$, the number $\|\vv v_1\we\dots\we\vv v_r\|$ will denote the supremum norm of the exterior product $\vv v_1\we\dots\we\vv v_r$ with respect to the standard basis of $\bigwedge^r(\R^k)$. Up to sign the coordinates of $\vv v_1\we\dots\we\vv v_r$ can be computed as all the $r\times r$ minors of the matrix composed of the coordinates of $\vv v_1,\dots,\vv v_r$ in the standard basis. Recall that the norm $\|\vv v_1\we\dots\we\vv v_r\|$ is {equivalent} to the $r$-dimensional volume of the parallelepiped
spanned by $\vv v_1,\dots,\vv v_r$, which is precisely the Euclidean norm of $\vv v_1\we\dots\we\vv v_r$\footnote{{In Diophantine approximation, if $\vv v_1,\dots,\vv v_r$ is a basis of $\Z^k\cap V$, where $V={\rm Span}_\R(\vv v_1,\dots,\vv v_r)$, the Euclidean norm of $\vv v_1\we\dots\we\vv v_r$ is known as the {\em height} of the linear rational subspace $V$ of $\R^k$, see \cite{MR1176315}.}}.

\bigskip

\begin{theorem}[\textbf{Quantitative non-Divergence estimate \cite[Theprem~5.2]{MR1652916}}]\label{QnD:t}
Let $k,d\in\N$, {$C,\alpha>0$, $0 < \rho\le 1/k$},  {a ball $B$} in $\R^d$
and {a function} $h:3^{k}B \to \GL_k(\R)$ be given. Assume that for
any linearly independent collection of integer vectors $\vv a_1,\dots,\vv a_r\in \Z^{k}$\\[-3ex]
\begin{enumerate}
\item[{\rm(i)}] the function $x\mapsto \|h(x)\vv a_1\we\dots\we h(x)\vv a_r\|$ is
$(C,\alpha)$-good\ on $3^{k}B$, \\[-1ex]
\item[{\rm(ii)}] $\sup\limits_{x\in B}\|h(x)\vv a_1\we\dots\we h(x)\vv a_r\|\ge\rho$.
\end{enumerate}
\vspace*{0ex}

Then for any $ \ve >0$
\begin{align}
\lambda_d\big(\big\{x\in B:\delta\big(h(x)\Z^{k}\big)<
\ve\big\}\big)
&\hspace*{1ex}\le\hspace*{1ex} kC (3^dN_d)^k \left(\frac\ve \rho
\right)^\alpha \lambda_d(B)\,,\label{QnD}
\end{align}
where $N_d$ is the Besicovitch constant.
\end{theorem}

%\bigskip

\noindent {We note that in
%To be absolutely precise
%\cite{MR1652916} contains the extra assumption that $\rho<1/k$. However, it was shown in
\cite{MR2434296} the second-named author established a version of the QnD estimate
%established in
where the
norm $\|\cdot\|$ is made to be Euclidean, not supremum. This made it possible to remove the condition $\rho<1/k$ and at the same time
%that this constrain was redundant. Also in \cite{MR2434296} condition (ii) was
replace condition (ii) above by a weaker condition
\begin{enumerate}
\item[{\rm(iii)}]  $\sup\limits_{x\in B}\|h(x)\vv a_1\we\dots\we h(x)\vv a_r\|\ge\rho^r$.
\end{enumerate}
The latter has been especially useful for studying \da\ on affine subspaces. See \cite{MR2648694} for a detailed exposition of the proof of the refined version of Theorem \ref{QnD:t} established in \cite{MR2434296}.}%and  \cite{MR3816502} for a further refinement.}

%\comm{Vitya, I don't think it's true. In \cite{MR2434296} I was working with Euclidean norm, hence the constant could be dropped. --- I have to think about it, I though it's ok regardless of the norm, of course the constant in \eqref{QnD} would change when you change the norm, but maybe I am missing something}

Estimate \eqref{QnD} is amazingly general. However, in many applications analysing conditions (i) and (ii){/(iii)} represents a substantial, often challenging, task. {However, in the case when $h$ is analytic, that is every entry of $h$ is a real analytic function of several variables,  condition (i) always holds for some $C$ and $\alpha$, see \cite{MR1652916} for details.} For the rest of this survey we shall mainly discuss developments in the theory of Diophantine approximation where the QnD estimate played an important, if not crucial, role. Naturally, we begin with the application of the QnD estimate that motivated  its discovery.

\section{The Baker-Sprind\v zuk conjecture and extremality}

In this section we explain the role of QnD in establishing the Baker-Sprind\v zuk conjecture -- a combination of two prominent problems in the theory of Diophantine approximation on manifolds, one due to A.~Baker \cite{MR0422171} and the other due to V.~G.~Sprind\v zuk \cite{MR586190}. The Baker-Sprind\v zuk conjecture is not merely a combination of two disjoint problems. Indeed, the origin of both problems lies in a single conjecture of Mahler \cite{MR1512754}  -- an important problem {in} the theory of transcendence posed in 1932 and proved by Sprind\v zuk in 1964 \cite{MR0245527}. Mahler's conjecture/Sprinz\v uk's theorem states that for any $n\in\N$ and any $\ve>0$ for almost every real number $x$ the inequality
\begin{equation}\label{Mah}
  |p+q_1x+q_2x^2+\dots+q_nx^n|<\|\vv q\|^{-n(1+\ve)}
\end{equation}
has only finitely many solutions  {$(p,\vv q)\in\Z\times\Z^{n}$, where $\vv q=(q_1,\dots,q_n)$ and}, as before, $\|\vv q\|=\max_{1\le i\le n}|q_i|$.
Note that, by Dirichlet's theorem or Minkowski's theorem for systems of linear forms, if $\ve\le0$ then \eqref{Mah} has infinitely many integer solutions $(p,\vv q)$ for any $x\in\R$. The condition $\ve>0$ in Mahler's conjecture is therefore sharp. It does not mean however that one cannot improve upon it! The improvements may come about when one replaces the right hand side of \eqref{Mah} with a different function of $\vv q$. One such improvement was conjectured by %Alan
Baker \cite{MR0422171} who proposed that the statement of Mahler's conjecture had to be true if the right hand side of \eqref{Mah} was replaced by
$$
\Pi_+(\vv q)^{-(1+\ve)},\qquad\text{where } \Pi_+(\vv q):=\prod_{\substack{i=1\\ q_i\neq0}}^n|q_i|\,.
$$
Clearly this leads to a stronger statement than Mahler's conjecture since
\begin{equation}\label{e202}
\Pi_+(\vv q)\le \|\vv q\|^n\,.
\end{equation}
The essence of replacement of $\|\vv q\|^n$ with $\Pi_+(\vv q)$ is to make the error of approximation depend on the size of each coordinate of $\vv q$ rather than on their maximum.

In another direction Sprind\v zuk \cite{MR586190} proposed a generalisation of Mahler's conjecture by replacing the powers of $x$ in \eqref{Mah} with arbitrary analytic functions of real variables which together with 1 are linearly independent over $\R$. The two conjectures (of Baker and Sprind\v zuk) can be merged in an obvious way to give what is known by the name of the Baker-Sprind\v zuk conjecture:

\medskip

\noindent{\textbf{The Baker-Sprind\v zuk Conjecture:}
{\em Let $f_1,\dots,f_n:U\to\R$ be {real analytic} functions
%of real variables
 defined on a connected open set $U\subset\R^d$. Suppose that $1,f_1,\dots,f_n$ are linearly independent over $\R$. Then for any $\ve>0$ {and} for almost every $x\in U$ the inequality
$$
  |p+q_1f_1(x)+\dots+q_nf_n(x)|<\Pi_+(\vv q)^{-1-\ve}
$$
has only finitely many solutions $(p,\vv q)\in\Z\times\Z^n$.
%, where $\vv q=(q_1,\dots,q_n)$.
}}

\medskip

At the time these conjectures were posed each seemed intractable and for a long while only limited partial results were known. Indeed, both conjectures remained open for $n\ge4$ almost until the turn of the millennium when they were solved in \cite{MR1652916} as an elegant application of the Quantitative non-Divergence estimate (Theorem~\ref{QnD:t}). The ultimate solution applies to the wider class of non-degenerate maps (defined below) that are not necessarily analytic.

\textbf{Non-degeneracy.} Let $\vv f=(f_1,\dots,f_n):U\to\R^n$ be a map defined on an open subset $U$ of $\R^d$.
Given a point $x_0\in U$, we say that $\vv f$ is \emph{$\ell$-non-degenerate at} $x_0$ if $\vv f$  is $\ell$ times continuously differentiable on some sufficiently small ball centered at $x_0$ and the partial derivatives
of $\vv f$ at $x_0$ of orders up to $\ell$ span $\R^n$. The map $\vv f$ is
called \emph{non-degenerate} at $x_0$ if it is $\ell$-non-degenerate at $x_0$ for some $\ell\in\N$; $\vv f$ is called non-degenerate almost everywhere (in $U$) if it is non-degenerate at almost every $x_0\in U$ with respect to Lebesgue measure. The non-degeneracy of differentiable submanifolds of $\R^n$ is defined via their parameterisation(s). Note that {a real  analytic} map $\vv f$ defined on a {connected open set} is non-degenerate {almost everywhere} if and only if $1,f_1,\dots,f_n$ are linearly independent over $\R$.

With the definition of non-degeneracy in place we are now ready to state the following flagship result of \cite{MR1652916} that solved the Baker-Sprind\v zuk conjecture in full generality not only in the analytic case, but also for arbitrary non-degenerate maps.

\medskip

\begin{theorem}[Theorem~A in \cite{MR1652916}]\label{KM98}
Let $\vv f=(f_1,\dots,f_n)$ be a map defined on an open subset $U$ of $\R^d$ which is non-degenerate almost everywhere. Then for any $\ve>0$, for almost every $x\in U$ {the inequality}
\begin{equation}\label{km1}
|p+q_1f_1(x)+\dots+q_nf_n(x)|< \Pi_+(\vv q)^{-1-\ve}
\end{equation}
has only finitely many solutions  {$(p,\vv q)\in\Z\times\Z^{n}$.}
\end{theorem}

%\bigskip
It is worth making further comments on the terminology used around the Baker-Sprind\v zuk conjecture. The point $\vv y\in\R^n$ is called {\em very well approximable (VWA)} if for some $\ve>0$ {the inequality}
\begin{equation}\label{km3}
|p+q_1y_1+\dots+q_ny_n|< %\Big(\max_{1\le i\le n}|q_i|\Big)
\|\vv q\|^{-n(1+\ve)}
\end{equation}
has infinitely many solutions $(p,\vv q)\in\Z^{n+1}%_{\neq0}
$.
The point $\vv y$ that is not VWA is often referred to as {\em extremal.}
The point $\vv y\in\R^n$ is called {\em very well multiplicatively approximable (VWMA)} if for some $\ve>0$  {the inequality}
\begin{equation}\label{km4}
|p+q_1y_1+\dots+q_ny_n|< \Pi_+(\vv q)^{-1-\ve}
\end{equation}
has infinitely many solutions $(p,\vv q)\in\Z^{n+1}$.
The point $\vv y$ that is not VWMA is often referred to as {\em strongly extremal.} As one would expect from a well set terminology we have that
$$
\vv y\text{ is strongly extremal}\qquad\Longrightarrow\qquad \vv y\text{ is extremal}\,.
$$
This is due to \eqref{e202}.

{Similarly, if $\mu$ is a measure on $\R^n$, one says that $\mu$ is   extremal or   strongly extremal if so is $\mu$-a.e.\ point of $\R^n$. The same goes for subsets of $\R^n$ carrying naturally defined measures. For example,  the notion of `almost all' for points lying on a differentiable submanifold in $\R^n$ can be defined in several  equivalent ways. Perhaps the simplest is to fix a parameterisation
$\vv f:U\to\R^n$ (possibly restricting $\cM$ to a local coordinate chart) and consider the pushforward
$\vv f*\lambda_d$ of Lebesgue measure on $\R^d$. Then a subset $S$ of $\vv f(U)$ is null {if and only if} $\lambda_d\big(\vv f^{-1}(S)\big)=0$.}
Now Theorem~\ref{KM98} can be rephrased as follows: {\em almost every point of any non-degenerate{\footnote{{Here we say that $\cM$ is non-degenerate if $\vv f$ is non-degenerate at $\lambda_d$-almost every point of $U$.}}} submanifold $\cM$ of $\R^n$ is strongly extremal}, or alternatively {\em any non-degenerate submanifold $\cM$ of $\R^n$ is strongly extremal}.

\bigskip

\noindent\textbf{Sketch of the proof of Theorem~\ref{KM98}} (for full details see \cite{MR1652916}). Define the following $(n+1)\times(n+1)$ matrix
\begin{equation}\label{new1}
u_{\vv f(x)}=\left(\begin{array}{ccc}
             1 & \vv f(x) \\[1ex]
             0 & I_n
           \end{array}
\right)
\end{equation}
and for $\vv t=(t_1,\dots,t_n)\in\Z^n_{\ge0}$ define the following $(n+1)\times(n+1)$ diagonal matrix
$$
g_{\vv t}=\left(\begin{array}{cccc}
            e^{t} &  &  & \\
             & e^{-t_1} &  & \\
              &  & \ddots &\\
              &  & & e^{-t_n}
          \end{array}
\right),\qquad \text{where }\ t=t_1+\dots+t_n.
$$
Given a solution $(p,\vv q)\in\Z^{n+1}$ to \eqref{km1}, one
defines $r=\Pi_+(\vv q)^{-\ve/(n+1)}$ and the smallest non-negative integers $t_i$ such that
$$
e^{-t_i}\max\{1,|q_i|\}\le r\,.
$$
Observe that $e^{t_i}<r^{-1}e\max\{1,|q_i|\}$. Then, by \eqref{km1},
$$
e^t|p+q_1f_1(x)+\dots+q_nf_n(x)|< {\big(e^nr^{-n}\Pi_+(\vv q)\big)}\Pi_+(\vv q)^{-1-\ve}=e^nr.
$$
Also, an elementary {computation} shows that $r<e^{-t\gamma}e^{n\gamma}$ with $\gamma=\ve/(n+1+n\ve)$. As a result we have that
\begin{equation}\label{km2}
\delta(g_{\vv t}u_{\vv f(x)}\Z^{n+1})<e^nr<e^{n(1+\gamma)}e^{-\gamma t}\,.
\end{equation}
Clearly, if for some $x\in U$ \eqref{km1} holds for infinitely many $\vv q$, then \eqref{km2} holds for infinitely many $\vv t\in{\Z^n_{\ge0}}$. Then the obvious line of proof of Theorem~\ref{KM98} is to demonstrate that the sum of measures
$$
\sum_{\vv t}\lambda_d\big\{x\in B:\delta(g_{\vv t}u_{\vv f(x)}\Z^{n+1})< e^{n(1+\gamma)} e^{-\gamma t}\big\}
$$
converges, where $B$ is a sufficiently small ball centred at an arbitrary point $x_0\in U$ such that $\vv f$ is non-degenerate at $x_0$. Indeed, the non-degeneracy condition placed on $\vv f$ justifies the restriction of $x$ to $B$ while the Borel-Cantelli Lemma from probability theory ensures that for almost all $x$ \eqref{km2} holds only finitely often subject to the convergence of the above series.

\vfil\eject

\noindent\underline{Remaining steps}:

$\bullet$ Take $h(x)=g_{\vv t}u_{\vv f(x)}$;

$\bullet$ Verify conditions (i) and (ii) of Theorem~\ref{QnD:t} for suitably chosen balls $B$;

$\bullet$ Conclude that
$$
\lambda_d(\{x\in B:\delta(g_{\vv t}u_{\vv f(x)}\Z^{n+1})< e^{n(1+\gamma)} e^{-\gamma t}\})\le \textsc{Const}\,\cdot\, e^{-\gamma \alpha t};
$$

$\bullet$ Observe that for each $t$ there are no more that $t^{n-1}$ integer $n$-tuples $\vv t$ such that $t_1+\dots+t_n=t$ and that
$$
\sum_{t=0}^\infty t^{n-1} e^{-\gamma \alpha t}<\infty
$$
thus completing the proof. \qed

\medskip

%\begin{remark}
It should be noted that the above sketch {of} proof is missing the details of verifying conditions (i) and (ii) of Theorem~\ref{QnD:t}. In particular, this requires explicit calculations of actions of $h$ on discrete subgroups of $\Z^{n+1}$ and understanding why certain maps are $(C,\alpha)$-good.
Details can be found in \cite{MR1652916}.
%\end{remark}

\subsection{Further remarks}\label{sec2.1}
The ideas of \cite{MR1652916} have been taken to a whole new level in \cite{MR2134453} by
%introducing the notion of (strongly) extremal measures in $\R^n$ and
identifying a large class of so-called {\em friendly measures} that are strongly extremal.
%Formally, a measure $\mu$ over $\R^n$ is called {\em (strongly) extremal} if the set of points $\vv x\in\R^n$ that are (strongly) extremal has $\mu$-measure zero; or in other words if $\mu$-almost every point in the support of $\mu$ is (strongly) extremal.
Examples include measures supported on a large class of fractal sets (the Cantor ternary set, the Sierpinski gasket, the attractors of certain IFSs (Iterated Function Systems), etc) and their pushforwards by non-degenerate maps. {See also    \cite{MR3816502} for further extension of the class of measures to which the QnD method applies.}

Theorem~\ref{KM98} was generalised in \cite{MR1982150} to affine subspaces of $\R^n$ (lines, hyperplanes, etc) that satisfy a natural Diophantine condition, as well as
%non-degenerate
to submanifolds of such affine subspaces. Affine subspaces and their submanifolds represent a very natural (if not the only natural) example of degenerate submanifolds of $\R^n$. One striking consequence of \cite{MR1982150} is the following criterion for analytic submanifolds: {\em an analytic submanifold $\cM$ of $\R^n$ is (strongly) extremal if and only if the {smallest} affine subspace of $\R^n$ that contains $\cM$ is (strongly) extremal.} {See also \cite{MR2434296} and \cite{dichotomy} for an extension of this `inheritance' principle to arbitrary (possibly non-extremal) affine subspaces  of $\R^n$.} Other natural generalisations of Theorem~\ref{KM98} include Diophantine approximation on complex analytic manifolds \cite{MR2094125}, Diophantine approximation in positive characteristics \cite{MR2321374}, and {in} $\Q_p$ ($p$-adic number{s}) and more generally Diophantine approximation {on submanifolds in} the product of several real and $p$-adic spaces \cite{MR2314053}.

The Baker-Sprind\v zuk conjecture deals with small values of one linear form of integer variables. {This is a} special case of the more general framework of systems of several linear forms in which the notions of extremal and strongly extremal matrices are readily available. Given an $n\times m$ matrix $X$ with real entries, one says that $X$ is extremal (not VWA) if and only if for any $\ve>0$
\begin{equation}\label{extremal}
\|\vv qX-\vv p\|^m< \|\vv q\|^{-n(1+\ve)}
\end{equation}
holds for at most finitely many $(\vv p,\vv q)\in\Z^m\times\Z^n$. The choice of the norm $\|\cdot\|$ does not affect the notion, but again for simplicity we choose the supremum norm: $\|\vv q\|=\max_{1\le i\le n}|q_i|$ for $\vv q=(q_1,\dots,q_n)$. Similarly, one says that $X$ is strongly extremal (not VWMA) if and only if for any $\ve>0$
\begin{equation}\label{extremal+}
\Pi(\vv qX-\vv p)< \Pi_+(\vv q)^{-(1+\ve)}
\end{equation}
holds for at most finitely many $(\vv p,\vv q)\in\Z^m\times\Z^n$. Here $\Pi(\vv y)=\prod_{j=1}^m|y_j|$ for $\vv y=(y_1,\dots,y_m)$ and, as before,
$\Pi_+(\vv q)=\prod_{i=1,~q_i\neq0}^n|q_i|$ for $\vv q=(q_1,\dots,q_n)$.

The QnD estimate can be applied to establish (strong) extremality of submanifolds of {the space of} matrices;
%, see \cite{MR2679461}. H
however, conditions (i) and (ii) of Theorem~\ref{QnD:t} are more difficult to translate into a `natural' and `practically checkable' definition of non-degeneracy. Indeed, identifying natural generalisations of the notion of non-degeneracy for matrices has been an active area of recent research, see \cite{MR2679461, MR3346961, Simmons1}. In the case of analytic {submanifolds of the space of} matrices  the goal has been attained %for extremal matrices
{in \cite{MR3777412}, with the notion of `constraining pencils' in the space of matrices replacing affine hyperplanes of $\R^n$.}
%however the more involved case of strongly extremal matrices is far from completion, see \cite{MR2679461}, \cite{MR3346961}.

\section{Khintchine-Groshev type results}

The theory of extremality %we have
discussed in the previous section deals with Diophantine inequalities with the right hand side written as the function
\begin{equation}\label{extr}
  \psi_\ve(h)=h^{-1-\ve}
\end{equation}
of either $\Pi_+(\vv q)$ or $\|\vv q\|^n$,
%=\max_{1\le i\le n}|q_i|^n$,
see \eqref{extremal} and \eqref{extremal+}. Of course, there are other choices for the `height' function of $\vv q$. For instance, in the case of the so-called {\em weighted} Diophantine approximation one uses the quasi-norm defined by
\begin{equation}\label{quasi}
\|\vv q\|_{\vv r}=\max_{1\le i\le n}|q_i|^{1/r_i}\,,
\end{equation}
where $\vv r=(r_1,\dots,r_n)\in\R^n_{>0}$ is an $n$-tuple of `weights' that satisfy the condition
\begin{equation}\label{weights}
  r_1+\dots+r_n=1\,.
\end{equation}
It this section we discuss the refinement of the theory of extremality that involves replacing the specific function $\psi_\ve$ given by \eqref{extr} with an arbitrary (monotonic) function $\psi$, akin to classical results of Khintchine \cite{MR1512207,MR1544787}. Below we state Khintchine's theorem in the one-dimensional case. Given a function $\psi:\N\to\R_+$, let
  $$
  \cA(\psi): =
\{x\in [0,1] \colon |qx-p|<\psi(q) \text{ for infinitely many }
(p,q)\in\Z\times\N \}  .
$$

\begin{theorem}[Khintchine's theorem]\label{Khi}
$$
    \lambda_1\big(\cA(\psi)\big) \ = \ \begin{cases} 0
      &\text{if } \quad \sum_{h=1}^\infty \, \psi(h)<\infty \ , \\[1ex]
      1
      &\text{if } \quad  \sum_{h=1}^\infty \, \psi(h)=\infty
      \  \text{ and %$q\psi(q)$
      {$\psi$ is
                    non-increasing}}.
                  \end{cases}
$$
\end{theorem}

%\comm{Vitya, are you sure? usually the assumption was that $\psi$ is monotonic, wasn't it?}

This beautiful finding has been generalised in many ways, and the theory for independent variables is now in a very advanced state, see for instance \cite{MR2508636} and \cite{MR2576284}. The generalisation of Khintchine's theorem to systems of linear forms was first established by Groshev \cite{Groshev-1938}. In the modern days theory, various generalisations of Khintchine's theorem to Diophantine approximation on manifolds are often called Khintchine type or Groshev type or Khintchine-Groshev type results. We will not define precise meaning of these words as there is some inconsistency in their use across the literature, although some good effort to harmonise the terminology was made in the monograph \cite{MR1727177}.

It has to be noted that the convergence case of Theorem~\ref{Khi} is a relatively simple consequence of the Borel-Cantelli Lemma. In the case of Diophantine approximation on manifolds this is no longer the case and establishing convergence Khintchine-Groshev type results for manifolds leads to a major challenge. Indeed, this was the case even in the special case associated with extremality. In this section we shall describe the role played by the QnD estimate in addressing this major challenge, namely in establishing convergence Khintchine-Groshev type refinements of Theorem~\ref{KM98} that were proved in \cite{MR1829381}. The key result of \cite{MR1829381} reads as follows.

\begin{theorem}[See \cite{MR1829381}]\label{KG1}
Let $\vv f=(f_1,\dots,f_n)$ be a map defined on an open subset $U$ of $\R^d$ which is non-degenerate almost everywhere. Let $\Psi:\Z^n\to\R_+$ be any function such that
\begin{equation}\label{cond1}
\Psi(q_1,\dots,q_i,\dots,q_n)\le\Psi(q_1,\dots,q'_i,\dots,q_n)\text{ if }|q_i|\ge |q'_i|\text{ and }q_iq_i'>0\,.
\end{equation}
Suppose that
\begin{equation}\label{conv1}
  \sum_{\vv q\in\Z^n}\Psi(\vv q)<\infty\,.
\end{equation}
Then for almost every $x\in U$ {the inequality}
\begin{equation}\label{km1p}
|p+q_1f_1(x)+\dots+q_nf_n(x)|< \Psi(\vv q)
\end{equation}
has only finitely many solutions $(p,\vv q)\in\Z^{n+1}$.
\end{theorem}

Observe that $\Psi(\vv q)=\Pi_+(\vv q)^{-1-\ve}$ {for} $\ve>0$ satisfies the conditions of Theorem~\ref{KG1}, and thus Theorem~\ref{KG1} is a true generalisation of Theorem~\ref{KM98}.
Prior to describing the ideas of the proof of Theorem~\ref{KG1} we formally state the following three corollaries: standard, weighted and multiplicative Khintchine-Groshev type results.

\begin{corollary}\label{KG1Cor1}
Let $\vv f$ be as in Theorem~\ref{KG1} and {let} $\psi:\R_+\to\R_+$ be any monotonic function. Suppose that
\begin{equation}\label{e022}
  \sum_{h=1}^\infty \psi(h)<\infty\,.
\end{equation}
Then for almost every $x\in U$ {the inequality}
\begin{equation}\label{e023}
|p+q_1f_1(x)+\dots+q_nf_n(x)|< \psi(\|\vv q\|^n)
\end{equation}
has only finitely many solutions $(p,\vv q)\in\Z^{n+1}$.
\end{corollary}

This Khintchine-Groshev type theorem, a direct generalisation of Sprind\v zuk's conjecture discussed in the previous section, is in fact a partial case of the following more general weighted version.

\bigskip

\begin{corollary}\label{KG1Cor2}
Let $\vv f$ be as in Theorem~\ref{KG1}, $\psi:\R_+\to\R_+$ be any monotonic function and
$\vv r=(r_1,\dots,r_n)\in\R^n_{>0}$ be an $n$-tuple satisfying \eqref{weights}. Suppose that
\eqref{e022} holds. Then for almost every $x\in U$ {the inequality}
\begin{equation}\label{e024}
|p+q_1f_1(x)+\dots+q_nf_n(x)|< \psi(\|\vv q\|_{\vv r})
\end{equation}
has only finitely many solutions $(p,\vv q)\in\Z^{n+1}$.
\end{corollary}

Recall again that Corollary~\ref{KG1Cor1} is a special case of Corollary~\ref{KG1Cor2}. Indeed, all one has to do to see it is to set $\vv r=(\frac1n,\dots,\frac1n)$. Note that Corollary~\ref{KG1Cor1} was established in \cite{MR1905790} using an approach that does not rely on the QnD estimate. However, without new ideas that approach does not seem to be possible to extend to the weighted case, let alone multiplicative approximation (the next corollary), where the QnD estimate has proven to be robust.

\begin{corollary}\label{KG1Cor3}
Let $\vv f$ be as in Theorem~\ref{KG1} and $\psi:\R_+\to\R_+$ be any monotonic function. Suppose that
\begin{equation}\label{e022+}
  \sum_{h=1}^\infty (\log h)^{n-1}\psi(h)<\infty\,.
\end{equation}
Then for almost every $x\in U$ {the inequality}
\begin{equation}\label{e023+}
|p+q_1f_1(x)+\dots+q_nf_n(x)|< \psi\big(\Pi_+(\vv q)\big)
\end{equation}
has only finitely many solutions $(p,\vv q)\in\Z^{n+1}$.
\end{corollary}

\bigskip

\noindent\textbf{Sketch of the proof of Theorem~\ref{KG1}} (for full details see \cite{MR1829381}). The proof again uses the QnD estimate (or rather an appropriate generalisation of Theorem~\ref{QnD:t}). However, this time the QnD estimates are not directly applicable to get the required result. The reason for that is that $\alpha$ in \eqref{QnD} is not matching the heuristic expectation, in fact it can hardly match it. The idea of the proof, which goes back to Bernik's paper \cite{MR1045454}\footnote{{It is worth mentioinng that in \cite{MR1045454} Bernik essentially proved Corollary~\ref{KG1Cor1} in the case $f_i(x)=x^i$.}}, is to separate two independent cases as described below.

\noindent\textbf{Case I.} Fix a small $\delta>0$ and consider the set of $x\in U$ such that \eqref{km1p} is satisfied simultaneously with the following condition on the gradient
\begin{equation}\label{km1+}
\big\|\nabla\big(q_1f_1(x)+\dots+q_nf_n(x)\big)\big\|\ge \|\vv q\|^{0.5+\delta}
\end{equation}
for infinitely many $(p,\vv q)\in\Z^{n+1}$.

\medskip

\noindent\textbf{Case II.} Fix a small $\delta>0$ and consider the set of $x\in U$ such that \eqref{km1p} is satisfied simultaneously with the opposite condition on the gradient:
\begin{equation}\label{km1+2}
\big\|\nabla\big(q_1f_1(x)+\dots+q_nf_n(x)\big)\big\|< \|\vv q\|^{0.5+\delta}
\end{equation}
for infinitely many $(p,\vv q)\in\Z^{n+1}$.

\medskip

Clearly, Theorem~\ref{KG1} would follow if one could show that the set $x\in U$ under consideration in each of the two cases is null. So how does that splitting into two cases help?

\medskip

Regarding Case II, note that the two conditions \eqref{cond1} and \eqref{conv1} imposed on $\Psi$ imply that
$$
\Psi(\vv q)\le \Pi_+(\vv q)^{-1}
$$
when $\|\vv q\|$ is sufficiently large. Hence it suffices to show that the set of $x\in U$ such that
\begin{equation}\label{km1+3}
\left\{\begin{array}{l}
|p+q_1f_1(x)+\dots+q_nf_n(x)|< \Pi_+(\vv q)^{-1}\\[1ex]
\big\|\nabla\big(q_1f_1(x)+\dots+q_nf_n(x)\big)\big\|< \|\vv q\|^{0.5+\delta}
       \end{array}\right.
\end{equation}
for infinitely many $(p,\vv q)\in\Z^{n+1}$ is null. For $\delta<0.5$ this effectively brings us back to an extremality problem, this time for matrices, which is dealt with using the QnD estimate in a similar way to the proof of Theorem~\ref{KM98}, albeit with much greater technical difficulties in verifying conditions (i) and (ii) of Theorem~\ref{QnD:t} (or rather an appropriate generalisation of Theorem~\ref{QnD:t}).

\medskip

Regarding Case I, the presence of the extra condition, namely condition \eqref{km1+}, leads to the following two key properties, which can be verified provided $U$ is of a sufficiently small (fixed) diameter.

\medskip

\noindent\textbf{Separation Property:} There is a constant $c>0$ such that for any $\vv q\in\Z^n$ with sufficiently large $\|\vv q\|$  for any $x\in U$ satisfying \eqref{km1+} the inequality
\begin{equation}\label{km+4}
|p+q_1f_1(x)+\dots+q_nf_n(x)|< \frac{c}{\big\|\nabla\big(q_1f_1(x)+\dots+q_nf_n(x)\big)\big\|}
\end{equation}
can hold for at most one integer value of $p$.

\medskip

\noindent\textbf{Measure Comparison Property:} For some constant $C>0$ for any $\vv q\in\Z^n$ with sufficiently large $\|\vv q\|$ and any integer $p$
\begin{equation}\label{km+6}
\lambda_d\big(\{x\in U: \text{\eqref{km1p} \&{} \eqref{km1+} hold}\}\big) \le C  \Psi(\vv q)   \lambda_d\big(\{x\in U: \text{\eqref{km1p} \&{} \eqref{km+4} hold}\}\big).
\end{equation}

\medskip
\medskip

The Separation Property implies that for a fixed $\vv q$
$$
\sum_{p\in\Z}\lambda_d\big(\{x\in U: \text{\eqref{km1p} \&{} \eqref{km+4} hold}\}\big)\le \lambda_d(U)\,.
$$
Putting this together with \eqref{km+6} and summing over $\vv q$ gives
$$
{\sum_{\substack{(p,\vv q)\in\Z^{n+1}\\ \vv q\neq\vv0}}}\lambda_d\big(\{x\in U: \text{\eqref{km1p} \&{} \eqref{km1+} hold}\}\big) \le C\lambda_d(U)\sum_{\vv q\in\Z^n}\Psi(\vv q)<\infty\,.
$$
It remain to apply the Borel-Cantelli Lemma to complete the proof.

\medskip

\begin{remark}
The above sketch proof is a significantly simplified version of the full proof presented in \cite{MR1829381}, which is far more effective. The effective elements of the proof in \cite{MR1829381} are stated as two independent results -- Theorems 1.3 and 1.4 in \cite{MR1829381}. These underpin a range of further interesting applications of {the} QnD {estimate} which we shall touch upon in later sections.
\end{remark}

\subsection{Further remarks}

Similarly to the theory of extremality, Khintchine-Groshev type results are not limited to non-degenerate manifolds and have been extended to affine subspaces of $\R^n$, see for instance \cite{MR2156655, MR2669718, MR2812651}. Khintchine-Groshev type results in $p$-adic and more generally $S$-arithmetic setting received their attention too, see for instance \cite{MR2570319, MR2897731}. Another remarkable application of the QnD estimate initially discovered in \cite{MR2110504} for polynomials and then extended in \cite{MR2481999} to non-degenerate curves enables one to remove the monotonicity constrain on $\psi$ from Corollary~\ref{KG1Cor1}.

The state-of-the-art for Diophantine approximation of matrices is far less satisfactory, where we do not have sufficiently general convergence Khintchine-Groshev type results. Partial results include for instance matrices with independent columns, see \cite{BBB}. The key difficulty likely lies within the case of simultaneous Diophantine approximation, which boils down to counting rational points near manifolds and will be discussed in \S\ref{RP}.

Theorem~\ref{KG1} implies that, under the convergence assumption \eqref{conv1}, for almost every $x\in U$ there exists a constant $\kappa>0$ such that
\begin{equation}\label{km1p+}
|p+q_1f_1(x)+\dots+q_nf_n(x)|\ge\kappa \Psi(\vv q)
\end{equation}
for all $(p,\vv q)\in\Z^{n+1}$ with $\vv q\neq\vv0$. Clearly, the constant $\kappa$ depends on $x$. Let $\cB(\vv f,\Psi;\kappa)$ denote the set of $x\in U$ satisfying the above condition for the same $\kappa$. Thus, $\cB(\vv f,\Psi;\kappa)$ is the set of all $x\in U$ such that \eqref{km1p+} holds for all $(p,\vv q)\in\Z^{n+1}$ with $\vv q\neq\vv0$. The conclusion of Theorem~\ref{KG1} exactly means that $\lambda_d(U\setminus\cB(\vv f,\Psi;\kappa))\to0$ as $\kappa\to0^+$. In general, the measure of $U\setminus\cB(\vv f,\Psi;\kappa)$ is positive. It is therefore of interest to understand how fast it is converging to zero. Motivated by applications in network information theory, this question was explicitely posed in \cite{J2010}. The answer was provided in \cite[Theorem~3]{MR3545930} and reads as follows:
$$
\lambda_d(U\setminus\cB(\vv f,\Psi;\kappa))\le \delta\lambda_d(U)
$$
for any
$$
\kappa\le\min\left\{\kappa_0,\ C_0\delta\left(\sum_{\vv q\in\Z^n\setminus\{\vv 0\}}\Psi(\vv q)\right)^{-1},\ C_1\delta^{d(n+1)(2l-1)}\right\}\,,
$$
where $l$ is a parameter characterising the non-degeneracy of $\vv f$ and $\kappa,C_0,C_1$ are positive (explicitely computable) constants that depend only on $\vv f$. More recently this `quantitative version' of the Khinthcine-Groshev theorem for manifolds has been generalised to affine subspaces \cite{MR3980277}.

Given that Khintchine's theorem (Theorem~\ref{Khi}) treats both the case of convergence and divergence, the question of establishing divergence counterparts to the convergence statements of this section is very natural. In this respect we now state the following known result.

\begin{theorem}\label{KG2}
Let $\vv f$ be as in Theorem~\ref{KG1}, $\psi:\R_+\to\R_+$ be any monotonic function and
$\vv r=(r_1,\dots,r_n)\in\R^n_{>0}$ be an $n$-tuple satisfying \eqref{weights}. Suppose that
\begin{equation}\label{e025}
  \sum_{h=1}^\infty \psi(h)=\infty\,.
\end{equation}
Then for a.e.\ $x\in U$ {the inequality \eqref{e024}}
has infinitely many solutions $(p,\vv q)\in\Z^{n+1}$.
\end{theorem}

Motivated by a Khintchine-Groshev generalisation of Mahler's conjecture, Theorem~\ref{KG2} was first established in \cite{MR1709049} for equal weights $r_1=\dots=r_n$ and functions $f_i(x)=x^i$ of one variable. Subsequently the method was generalised in \cite{MR1944505} to arbitrary non-degenerate maps
and more lately to arbitrary weights \cite{MR2989975}.
Theorem~\ref{KG2} provides the divergence counterpart to Corollaries~\ref{KG1Cor1} and \ref{KG1Cor2}. Establishing the divergence counterpart to Corollary~\ref{KG1Cor3} (the multiplicative case) remains a challenging open problem even in dimension $n=2$:

\medskip

\begin{problem}\label{p1}
Let $\vv f$ be as in Theorem~\ref{KG1} and $\psi:\R_+\to\R_+$ be any monotonic function. Suppose that
\begin{equation}\label{e022++}
  \sum_{h=1}^\infty (\log h)^{n-1}\psi(h)=\infty\,.
\end{equation}
Prove that for almost every $x\in U$ \eqref{e023+} holds for infinitely many $(p,\vv q)\in\Z^{n+1}$.
\end{problem}

In the light of the recent progress on a version of Problem~\ref{p1} for simultaneous rational approximations for lines made in \cite{BHV,SamChow,SamLei,SamNiclas} it would also be very interesting to investigate Problem~\ref{p1} when $\vv f$ is a linear/affine (and hence degenerate) map from $\R^d$ into $\R^n$ ($1\le d<n$). Note that the case $d=n$ of Problem~\ref{p1} follows from a result of Schmidt \cite{Schm}.

We conclude this section by one final comment: the proof of Theorem~\ref{KG2} is also underpinned by the QnD estimate. At first glance this may seem rather counter-intuitive since the QnD estimate deals with upper bounds, while Theorem~\ref{KG2} is all about lower bounds. One way or another, this is the case and we shall return to explaining the role of the QnD estimate in establishing Theorem~\ref{KG2} in the next section within the more general context of Hausdorff measures.

\section{Hausdorff measures and dimension}

It this section we discuss another refinement of the theory of extremality that aims at understanding the Hausdorff dimension of exceptional sets. The basic question is as follows: {\em given a non-degenerate submanifold $\cM$ of $\R^n$ and $\ve>0$ (not necessarily small), what is the Hausdorff dimension of the set of $\vv y\in\cM$ such that \eqref{km3} holds for infinitely many $(p,\vv q)\in\Z^{n+1}$.} The same question can be posed in the multiplicative setting \eqref{km4}, for weighted Diophantine approximation (when $\|\vv q\|^n$ is replace by $\|\vv q\|_{\vv r}$ given by \eqref{quasi}) and for Diophantine approximation of matrices such as \eqref{extremal} and \eqref{extremal+}.

The background to this question lies with the classical results of Jarn\'ik \cite{Ja29} and Besicovitch \cite{MR1574327} stated below in the one-dimensional case.

\begin{theorem}[The Jarn\'ik-Besicovitch theorem]\label{JB}
Let $\ve>0$ and $\psi_\ve(x)=x^{-1-\ve}$. Then
$$
  \dim \cA(\psi_\ve) =
       \ \frac2{2+\ve} \,.
$$
\end{theorem}

This %this
fundamental result tells us exactly how the size of $\cA(\psi_\ve)$ get smaller as we increase $\ve$ and thus make the approximation function $\psi_\ve$ decrease faster. As with Khintchine's theorem, the Jarn\'ik-Besicovitch theorem has been generalised in many ways and the theory for independent variables is now in a very advanced state, see for instance \cite{MR2259250}, \cite{MR2508636}, \cite{MR2576284} and \cite{MR3798593}. It has to be noted that showing the upper bound $\dim \cA(\psi_\ve)\le \frac2{2+\ve}$ is a relatively simple consequence of the so-called Hausdorff-Cantelli Lemma \cite{MR1727177}, an analogue of the Borel-Cantelli Lemma. However, in the case of Diophantine approximation on manifolds establishing upper bounds for manifolds leads to a major challenge that is very much open till these days (see Problem~\ref{p2} below). On the contrary, the lower bounds have been obtained in reasonable generality. The main purpose of this section is to exhibit the role played by the QnD estimate in obtaining the lower bounds.

As before $\vv f:U\to \R^n$, and $\vv r=(r_1,\dots,r_n)$ is an $n$-tuple of positive numbers satisfying \eqref{weights}. Given $s>0$, $\cH^s$ will denote the $s$-dimensional Hausdorff measure. Let $\delta>0$, $H>1$ and
$$
\Phi_{\vv r}(H,\delta):=\left\{x\in U: \exists\ \vv q\in\Z^n\smallsetminus\{0\} \text{ such that }
\left\{\begin{array}{l}
|p+\vv q\cdot\vv f(x)|<\delta H^{-1}\\
\|\vv q\|_{\vv r}\le H
       \end{array}
\right.\right\}.
$$
Also, given a function $\psi:\R_+\to\R_+$, let $\cW(\vv f,U,\vv r,\psi)$ be the set of $x\in U$ such that
\eqref{e024} has infinitely many solutions $(p,\vv q)\in\Z^{n+1}$. The following homogeneous version of one of the main results from \cite{MR2989975} works as a black-box to proving divergence Khintchine-Groshev type results (such as Theorem~\ref{KG2}) and lower bounds for Hausdorff dimension.

\begin{theorem}[Theorem~3 in \cite{MR2989975}]\label{tnice}
Let $\vv f:U\to\R^n$ be a $C^2$ map on an open subset $U$ of $\R^d$, and {let} $\vv r=(r_1,\dots,r_n)$ be an $n$-tuple of positive numbers satisfying \eqref{weights}. Suppose that for almost every point $x_0\in U$ there is an open neighborhood $V\subset U$ of $x_0$ and constants $0<\delta,\omega<1$ such that for any ball $B\subset V$ we have that
\begin{equation}\label{nice}
\lambda_d\big(\Phi_{\vv r}(H,\delta)\cap
B\big)\le \omega\lambda_d(B)
\end{equation}
for all sufficiently large $H$.
Let $d-1<s\le d$ and $\psi:\R_+\to\R_+$ be monotonic. Then
$$
\cH^s\big(\cW(\vv f,U,\vv r,\psi)\big)=\cH^s(U)\qquad\mbox{if }
\qquad \sum_{\vv q\in\Z^n_{\neq0}} \|\vv q\|  \left(\frac{\psi(\|\vv q\|_{\vv r})}{\|\vv q\| }\right)^{s+1-d}=\infty.
$$
\end{theorem}

The proof of this result makes use of {\em ubiquitous systems} as defined in \cite{MR2184760}, see also \cite{MR1773667}, \cite{MR1944505} and the survey \cite{MR1975457} for the related notion of {\em regular systems}. The QnD estimate steps in when one wishes to apply Theorem~\ref{tnice}, namely to verify condition \eqref{nice}. In particular, as was demonstrated in \cite{MR2989975}, any non-degenerate map $\vv f$ safisfies this condition for any collection of weights $\vv r$. Upon taking $s=d$ one then verifies that
$$
\text{Theorem~\ref{tnice}}\qquad\Longrightarrow\qquad\text{Theorem~\ref{KG2}}.
$$

Another consequence of Theorem~\ref{tnice} is the following lower bound on the Hausdorff dimension of exceptional sets that contributes to resolving the problem outlined at the beginning of this section. For simplicity we only state the result for the case of equal weights:  $\vv r_0=(\frac1n,\dots,\frac1n)$.

\begin{corollary}[Corollary~2 in \cite{MR2989975}]\label{cor1}
Let $\vv f$ be as in Theorem~\ref{KG1}, $\psi:\R_+\to\R_+$ be any monotonic function and let
$$
\tau_{\psi}:=\liminf_{t\to\infty}\frac{-\log\psi(t)}{\log t}\,.
$$
Suppose that
$n  \leq \tau_\psi  < \infty $. Then
\begin{equation}\label{dimdim}
\dim \cW(\vv f,U,\vv r_0,\psi)\ \ge \ s:=d-1+ \frac{n+1}{n\tau_\psi+1} \, .
\end{equation}
\end{corollary}

The number $\tau_{\psi}$ is often referred to as the {\em lower order} of $1/\psi$ at infinity. It indicates the growth of the function $1/\psi$ `near' infinity. Naturally for $\psi_\ve(t)=t^{-1-\ve}$ we have that $\tau_{\psi_\ve}=1+\ve$.
Estimate \eqref{dimdim} was previously shown in \cite{MR1738177} {for arbitrary extremal submanifolds of $\R^n$}. The additional benefit of Theorem~\ref{tnice} compared to \eqref{dimdim} is that it allows one to compute the Hausdorff measure of $\cW(\vv f,U,\vv r_0,\psi)$ at $s=\dim\cW(\vv f,U,\vv r_0,\psi)$. It is believed that the lower bound given by \eqref{dimdim} is exact for non-degenerate maps $\vv f$ at least in the analytic case. It is readily seen that to establish the desired equality for all $\psi$ in question it suffices to consider approximation functions $\psi_\ve$ only. Hence the following

\bigskip

\begin{problem}\label{p2}
Let $\vv f:U\to\R^n$ be an analytic non-degenerate maps defined on a ball $U$ in $\R^d$. Then for every $\ve>0$
\begin{equation}\label{dimdim+}
\dim \cW(\vv f,U,\vv r_0,\psi_\ve)\ = \ d-1+ \frac{n+1}{n+1+n\ve} \, .
\end{equation}
\end{problem}

\bigskip

Problem~\ref{p2} was established in full in the case $n=2$ by R.\,C.\ Baker \cite{Baker-1978}. For $n\ge3$ it remains very much open. However, for the polynomial maps $\vv f=(x,\dots,x^n)$ it was settled by Bernik in \cite{Bernik-1983a} for arbitrary $n$.

\subsection{Further remarks}\label{HDFR}

Although establishing the upper bound in \eqref{dimdim+} remains a prominent open problem, it was shown in \cite{MR2069553} that the QnD estimate can be used to resolve it for small $\ve$, namely for
$0<\ve <1 /(4n^2+2n-4)$ {when $d=1$}. This seemingly inconsequential result turned out to have major significance in resolving open problems on badly approximable points on manifolds, which will be discussed in Section \ref{BAP}.

Problem~\ref{p2} can be refined in the spirit of Khintchine-Groshev type result by asking to prove the following convergence counterpart to Theorem~\ref{tnice}.

\medskip

\begin{problem}\label{p3}
Let $\vv f:U\to\R^n$ be a non-degenerate analytic map defined on a ball $U$ in $\R^d$, $0<s<d$ and $\vv r=(r_1,\dots,r_n)$ be an $n$-tuple of positive numbers satisfying \eqref{weights}. Suppose that
$\psi:\R_+\to\R_+$ is any monotonic function.
Prove that
$$
\cH^s\big(\cW(\vv f,U,\vv r,\psi)\big)=0\qquad\mbox{if }
\qquad \sum_{\vv q\in\Z^n_{\neq0}} \|\vv q\|  \left(\frac{\psi(\|\vv q\|_{\vv r})}{\|\vv q\| }\right)^{s+1-d}<\infty.
$$
\end{problem}

For $n=2$ Problem~\ref{p3} was resolved in \cite{MR3824783} for equal weights but in higher dimensions {it is still open} even for $\vv f(x)=(x,\dots,x^n)$, let alone non-degenerate maps. In fact, for $n=2$ the case of non-equal weights can be reduced to the case of equal weights. This can be shown by modifying the proof of Theorem~2 from \cite{BB00}.

\section{Rational points near manifolds}\label{RP}

Rational and integral points lying on or near curves and surfaces crop up in numerous problems in number theory and are often one of the principle objects of study (e.g. in analytic number theory, Diophantine approximation, Diophantine geometry).  The goal of this section is to demonstrate the role of the QnD estimate in recent counting results on rational points lying close to manifolds  \cite{MR2874641, MR2373145, MR2729002, divergence}. The motivation lies within the theory of simultaneous Diophantine on manifolds, which boils down to understanding the proximity of rational points $\vv p/q=(p_1/q,\dots,p_n/q)$ to points
$\vv y=(y_1,\dots,y_n)\in \R^n$ restricted to a submanifold $\cM$. Here
$p_1,\dots,p_n\in\Z$ and $q\in\N$ is the common denominator of the coordinates of $\vv p/q$.
By Dirichlet's theorem, for every irrational point $\vv y\in\R^n$ there are infinitely many $\vv p/q\in\Q^n$ such that
$$
\max_{1\le i\le n}\left|y_i-\frac{p_i}q\right|<q^{-1-1/n}\,.
$$
This inequality can be rewritten in the form
\begin{equation}\label{e32}
\max_{1\le i\le n}\left|qy_i-p_i\right|^n<q^{-1}\,.
\end{equation}
Understanding when the right hand side of \eqref{e32} can be replaced by $\psi_\ve(q)=q^{-1-\ve}$, $\ve>0$ (or even by a generic monotonic function $\psi$) is the subject matter of many classical problems and famous results. The basic question is about the solvability of
\begin{equation}\label{e33}
\max_{1\le i\le n}\left|qy_i-p_i\right|^n<\psi_\ve(q)
\end{equation}
in $(\vv p,q)\in\Z^n\times\N$ for arbitrarily large $q$.
For example, celebrated Schmidt's subspace theorem states that for any algebraic $y_1,\dots,y_n$ such that $1,y_1,\dots,y_n$ are linearly independent over $\Q$ and any $\ve>0$ ~~\eqref{e33} has only finitely many solutions $(\vv p,q)\in\Z^n\times\N$.

When the point $\vv y$ lies on a submanifold $\cM$, \eqref{e33} forces the rational point $\vv p/q$ {to} lie near $\cM$. Hence, metric problems concerting \eqref{e33} (e.g. Khintchine type results, analogues of the Jarn\'ik-Besicovitch theorem, etc) have resulted in significant interest into counting and understanding the distribution of rational points lying close to submanifolds of $\R^n$. The basic setup is as follows. Given $Q>1$ and $0<\psi<1$, let
$$
R_\cM(Q,\psi)\ =\ \Big\{(\vv p,q)\in\Z^n\times\N:\ 1\le q\le Q,\ {\rm dist}(\vv
p/q,\cM)\le\psi/q\Big\}\,.
$$
It is not difficult to work out the following
\begin{equation}\label{e34}
\text{Heuristic:\qquad $\#R_\cM(Q,\psi) \asymp \psi^{m} Q^{d+1}$,}
\end{equation}
where $d=\dim\cM$, $m=n-d=\codim\cM$ and $\#$ stands for the cardinality. Also, $\asymp$ means the simultaneous validity of two Vinogradov symbols $\ll$ and $\gg$, {where} $\ll$ means the inequality $\le$ up to a positive multiplicative constant.

\medskip

\begin{remark}\label{r3}
This heuristic estimate has to be treated with caution as, for instance, the unit circle $y_1^2+y_2^2=1$ will always contain at least \textsc{const}$\cdot Q$ rational points (given by {Pythagorean} triples) resulting in $\#R_\cM(Q,\psi)\gg Q $ no matter how small $\psi$ is. On the contrary, the circle $y_1^2+y_2^2=3$ contains no rational points resulting in $\#R_\cM(Q,\psi)=0$ for large $Q$ when $\psi=o(Q^{-1})$. Also, rational (affine) subspaces inherently contain many rational points and so any manifold that contains a rational subspace may break the heuristic with ease for moderately small $\psi$.
\end{remark}

The following is the {principal} problem in this area, see \cite{MR2874641}, \cite{RP-Huang}.

\begin{problem}\label{p4}
Show that \eqref{e34} holds for any `suitably curved' compact differentiable submanifold $\cM$ of $\R^n$ when $\psi\ge Q^{-1/m+\delta}$, where $\delta>0$ is arbitrary and $m=\codim\cM$.
\end{problem}

Ideally, it would be desirable to resolve this problem for all non-degenerate submanifolds of $\R^n$. The condition $\psi\ge Q^{-1/m+\delta}$ is pretty much optimal unless one imposes further constraints on the internal geometry of $\cM$. For instance, to relax the condition on $\psi$ one has to exclude the manifolds that contain a rational subspace of dimension $d-1$ (when $d>1$). In what follows we shall describe the role of the QnD estimate in establishing the lower bound for analytic non-degenerate manifolds.

\begin{theorem}[Corollary~1.5 in \cite{MR2874641}]\label{RP:t}
For any analytic non-degenerate submanifold $\cM\subset\R^n$ of dimension $d$ and codimension $m=n-d$ there exist constants $C_1,C_2>0$ such that
\begin{equation}\label{e602}
\#R_\cM(Q,\psi) \ge  C_1 \psi^{m} Q^{d+1}
\end{equation}
for all sufficiently large $Q$ and all real $\psi$ satisfying
\begin{equation}\label{e600}
    C_2Q^{-1/m }< \psi < 1.
\end{equation}
\end{theorem}

\noindent\textbf{Sketch of the proof} (for full details see \cite{MR2874641}).
For simplicity we will assume that $\cM$ is the image $\vv f(\cU)$ of a cube $\cU\subset\R^d$ under a map $\vv f$ and that $\cM$ is bounded. Further, without loss of generality we will assume that $\vv f$ is of the {\em  Monge form}, that is
$$
\vv f(x_1,\dots,x_d)=\big(x_1,\dots,x_d,f_1(x_1,\dots,x_d),\dots,f_m(x_1,\dots,x_d)\big)\,.
$$
The rational points $\vv p/q$ give rises to the integer vectors $\vv a=(q,p_1,\dots,p_n)$, which are essentially projective representations of $\vv p/q$.
Define $\vv y(x)=\big(1,\vv f(x)\big)$, a projective representation of $\vv f(x)$. As is well known, the distance of $\vv p/q$ from $\vv f(x)$ is comparable to the projective distance between them defined as the sine of the acute angle between $\vv a$ and $\vv y(x)$. To make this angle $\ll\psi/q$ and thus ensure that $\vv p/q$ lies in $R_\cM(Q,\psi)$, it is enough to verify that
\begin{equation}\label{e37}
|q|\le Q\qquad\text{and}\qquad
|\vv g_i(x)\cdot\vv a|\ll \psi \qquad(1\le i\le n)
\end{equation}
for any fixed collection $\vv g_1(x),\dots,\vv g_n(x)$ of vector orthogonal to $\vv y(x)$ and such that $\|\vv g_i(x)\|\ll1$ and
$
\|\vv g_1(x)\we\dots\we\vv g_n(x)\|\asymp \|\vv g_1(x)\|\cdots\|\vv g_n(x)\|\,.
$
Clearly \eqref{e37} defines a convex body of $\vv a\in\R^{n+1}$ and one can potentially use Minkowski's theorem on convex bodies to find $\vv a$. However, for $\psi$ much smaller than $Q^{-1/n}$ this is impossible -- the volume of the body is too small. To overcome this difficulty, the convex body is expanded in the directions tangent to the manifolds written in the projective coordinates, see \eqref{e38} below. For this purpose it convenient to make the following choices for $\vv g_i(x)$. First of all,
for $i=1,\dots,m$ the vectors $\vv g_i(x)$ are taken to be orthogonal to $\vv y(x),\partial\vv y(x)/\partial x_1,\dots,\partial\vv y(x)/\partial x_d$, which are linearly independent for $\vv f$  of the Monge form.
Next, for $i=m+1,\dots,n$ the vectors $\vv g_i(x)$ are taken to be orthogonal to
$\vv y(x)$, $\vv g_1(x)$, \dots, $\vv g_{m}(x)$.

We will need two auxiliary positive parameters $\widetilde\psi$ and $\widetilde Q$ which will be proportional to $\psi$ and $Q$ respectively. Now, for $\kappa>0$ and a given $x$, we consider the convex body of $\vv a\in\V$ defined by
\begin{equation}\label{e38}
\begin{array}{l}
         |\vv g_i(x)\dt\vv a| < \widetilde\psi\quad(1\le i\le m)\,,\\[0.4ex]
         |\vv g_{m+j}(x)\dt\vv a| < (\widetilde\psi^{m } \widetilde Q)^{-\frac{1}{d}}\quad(1\le j\le d)\,,\\[0.7ex]
         |q|\le  \kappa \widetilde Q\,.
\end{array}
\end{equation}
For a suitably chosen constant $\kappa=\kappa_0$ dependent only on $n$ and $\vv f$, this body is of sufficient volume to apply Minkoswki's theorem. This results in the existence of an $\vv a\in\Z^{n+1}\smallsetminus\{0\}$ satisfying \eqref{e38}. Thus the set
$$
B(\widetilde\psi,\widetilde Q,\kappa)=\{x\in U: \exists\ \vv a\in\Z^{n+1}_{\neq0} \text{ satisfying \eqref{e38}}\}
$$
coincides with all of $U$ for $\kappa=\kappa_0$.
Now suppose that $\kappa_1<\kappa_0$ and
$$
x\in U\smallsetminus B(\widetilde \psi,\widetilde Q,\kappa_1)\qquad(~=B(\widetilde\psi,\widetilde Q,\kappa_0)\smallsetminus B(\widetilde\psi,\widetilde Q,\kappa_1)~).
$$
Then, the first two collections of inequalities in \eqref{e38} are satisfied for some {vector} $\vv a=(q,p_1,\dots,p_n)\in\Z^{n+1}$ such that
$$
\kappa_1 \widetilde Q\le |q|\le \kappa_0\widetilde Q\,.
$$
The first set of inequalities in \eqref{e38} keeps the rational point $\vv p/q$ at distance $\ll\ve_1=\widetilde \psi/\widetilde Q$ from the tangent plane to the manifold at $\vv f(x)$. The second set of inequalities in \eqref{e38} keeps the rational point $\vv p'/q$, where $\vv p'=(p_1,\dots,p_d)$, at distance $\ll \ve_2=(\widetilde \psi^{m } \widetilde Q)^{-\frac{1}{d}}\widetilde Q^{-1}$ from $x$. Since the tangent plane deviates from the manifold quadratically, assuming that $\ve_2^2\le\ve_1$ we conclude that the point $\vv p/q$ remains at distance $\ll \ve_1+\ve_2^2\ll\ve_1$ from the manifold, see \cite[Lemma~4.3]{MR2874641}. The proof of this uses nothing but {the second-order} Taylor's formula.
% to the second order.

To sum up,
$$
U\smallsetminus B(\widetilde \psi,\widetilde Q,\kappa_1)\subset\bigcup_{(\vv p,q)\in R_\cM(c_1\widetilde Q,c_2\widetilde \psi)} B(\vv p'/q,c_3\ve_2)\,,
$$
where $B(x,r)$ is a ball in $\R^d$ centered at $x$ of radius $r$, {and} $c_1,c_2,c_3>0$ {are} some constants. Then
\begin{equation}\label{e39}
\lambda_d\big(U\smallsetminus B(\widetilde \psi,\widetilde Q,\kappa_1)\big)\ll  \#R_\cM(c_1\widetilde Q,c_2\widetilde \psi)(c_3\ve_2)^d \asymp \#R_\cM(c_1\widetilde Q,c_2\widetilde \psi)\,\psi^{-m} Q^{-(d+1)}\,.
\end{equation}
At this point the QnD estimate is used to verify that for a suitably small constant $\kappa_1$ the set $B(\widetilde \psi,\widetilde Q,\kappa_1)$ has small measure, say, $\le\frac12\lambda_d(U)$. {Thus}
$\lambda_d\big(U\smallsetminus B(\widetilde \psi,\widetilde Q,\kappa_1)\big)\ge\tfrac12\mu_d(U)$, and \eqref{e39} implies the {desired} result on requiring that $\widetilde Q\le Q/c_1$ and $\widetilde \psi\le \psi/c_2$.

\medskip

To finish this discussion, we shall show explicitly how \eqref{e38} can be re-written for the purpose of applying the QnD {estimate}. For simplicity we consider a non-degenerate planar curve $\cC=\{(x,f(x):x\in U\}$, where $U$ is an interval, and so $d=m=1$ and $n=2$, and restrict ourselves to the case when $\widetilde \psi\asymp \widetilde Q^{-1}$. The latter means that we are counting rational points closest possible to $\cC$.
Let $\widetilde \psi \widetilde Q=\kappa^{1/3}$, $e^t=\widetilde \psi^{-1}\kappa^{1/3}$ and $e^{-t}=\widetilde Q^{-1}\kappa^{-2/3}$.

Then, \eqref{e38} can be replaced by
\begin{equation}\label{e41+}
\delta(g_tG_x\Z^3)<\kappa^{1/3},
\end{equation}
where
$$
G_x=\left(\begin{array}{ccc}
            f(x)-xf'(x) & f'(x) & -1 \\
            x & -1 & 0 \\
            1 & 0 & 0
          \end{array}
\right)
$$
and
$$
g_t=\left(\begin{array}{ccc}
            e^t & 0 & 0 \\
            0 & 1 & 0 \\
            0 & 0 & e^{-t}
          \end{array}
\right)\,.
$$
Indeed, the first row of $G_x$ is simply $\vv g_1(x)$ appearing in \eqref{e38}, and
the second row of $G_x$ is simply a multiple of $\vv g_2(x)$ appearing in \eqref{e38}.
Thus, counting rational points closets to a planar curve as discussed above relies upon finding an appropriately small constant $\kappa>0$ such that the set of $x\in U$ satisfying \eqref{e41+} has measure at most, say, $\tfrac12\lambda_1(U)$ for all sufficiently large $t$. To rephrase this, half of the curve $x\mapsto G_x\Z^3$ in $X_3$ has to remain in the compact set $\cK_{\ve}$ defined by \eqref{K_e} with $\ve=\kappa^{1/3}$ under the action by the $g_t$, that is
\begin{equation}\label{orbit}
\lambda_1(\{x\in U:g_tG_x\Z^3\in\cK_{\ve}\})~\ge~\tfrac12\lambda_1(U)\qquad\text{for all sufficiently large $t$}\,.
\end{equation}

\subsection{Further remarks}

When $d=1$ it was shown in \cite[Theorem~7.1]{MR2874641} that for analytic non-degenerate curves \eqref{e600} can be relaxed to
\begin{equation}\label{e601}
    C_2 Q^{-\frac{3}{2n-1}}< \psi < 1\,.
\end{equation}
More recently, the condition of the analyticity was removed in \cite{divergence} following a more careful and explicit application of the QnD estimate.  {In essence, the analytic case does not require to deal with condition (i) of Theorem~\ref{QnD:t}, the latter task being
%and dealing with this condition is the task
accomplished in \cite{divergence}.} For $d>1$ removing the analyticity condition from Theorem~\ref{RP:t} remains an open problem.
In the case of planar curves Problem~\ref{p4} was solved for non-degenerate planar curves as a result of \cite{MR1310631, MR2373145, MR2729002, MR2242634}, see also asymptotic and inhomogeneous results in \cite{MR2777039}, \cite{MR3318157}, \cite{MR3630723}, \cite{MR3263942}.
Upper bounds in higher dimensions represent a challenging open problem, but see \cite{MR3658127}, \cite{MR3809714} and \cite{RP-Huang} for some recent results.

For $n=2$ ($d=1$) condition \eqref{e601} does not actually improve upon \eqref{e600}. In fact, \eqref{e600} is optimal within the class of all non-degenerate hypersurfaces, and in particular non-degenerate planar curves, see Remark~\ref{r3} above. In principle, the existence of rational points as opposed to counting does not require using the QnD, see \cite{MR3731303}.

Detecting rational points near planar curves closer than the limit set by the left hand side of \eqref{e602} will require additional conditions on top of non-degeneracy and represents an interesting problem:

\begin{problem}\label{p5}
Find `reasonable' conditions on a connected analytic curve $\cC$ in $\R^2$ sufficient to satisfy
\begin{equation}\label{e41}
\liminf_{q\to\infty}q^2{\rm dist}\,(\cC,\tfrac1q\Z^2)=0\,.
\end{equation}
\end{problem}
Observe that for ellipses in $\R^2$ Problem~\ref{p5} reduces to {the Oppenheim} conjecture (1929) remarkably proved by Margulis in 1986:

\begin{theorem}[Margulis, 1986]\label{t9}
Let $Q$ be a nondegenerate indefinite quadratic form of 3 real variables and suppose that $Q$ is not a multiple of a form with rational coefficients. Then for any $\ve>0$ there exist nonzero integers $a,b,c$ such that
\begin{equation}\label{opp}
|Q(a,b,c)|<\ve.
\end{equation}
\end{theorem}

{To see the link between Theorem~\ref{t9} and Problem~\ref{p5}, first divide \eqref{opp} through by $c^2$ and, using the fact that $Q$ is a homogeneous polynomial of degree $2$,  obtain the following equivalent inequality:
\begin{equation}\label{opp1}
\Big|Q\Big(\frac ac,\frac bc,1\Big)\Big|<\frac{\ve}{c^{2}}.
\end{equation}
Since $Q$ is indefinite, without loss of generality one can assume that
$$Q(x,y,1)=\widetilde Q(x-x_0,y-y_0)-r^2$$ for some positive definite quadratic form $\widetilde Q$ of two variables and some $r>0$. If necessary one can permute the variables $a$, $b$ and $c$ to make sure this is the case. If $\cC$ denotes the curve in $\R^2$ defined by the equation $Q(x,y,1)=0$, then an elementary check shows that \eqref{opp1} is equivalent to ${\rm dist}\big(\cC,(a/c,b/c)\big)\ll\ve/c^2$, and also that $|c|\gg \max\{|a|,|b|\}$. Hence it becomes obvious that Theorem~\ref{t9} is equivalent to \eqref{e41} for the specific type of curves $\cC$ in question.
For instance, if $Q(x,y,z)=x^2+y^2-(rz)^2$ for some $r>0$, then $\cC$ is the circle of radius $r$ centred at the origin. In general, $\cC$ is an ellipse.}

%
%To see the link between Theorem~\ref{t9} and Problem~\ref{p5}
%one can assume, without loss of generality, that $c\ge\max\{|a|,|b|\}$. Indeed, this can be easily met by applying appropriate rotations of $\R^3$ by $90^\circ$ lying in $\operatorname{SO}(3,\Z)$ and the reflection with respect to the origin. Then, dividing \eqref{opp} through by $c^2$ and using the fact that $Q$ is a homogeneous polynomial of degree $2$ we obtain the following equivalent inequality
%\begin{equation}\label{opp1}
%\Big|Q\Big(\frac ac,\frac bc,1\Big)\Big|<\frac{\ve}{c^{2}}.
%\end{equation}
%If $\cC$ is the curve in $\R^2$ defined by the equation $Q(x,y,1)=0$, then
%an elementary check shows that \eqref{opp1} is equivalent to ${\rm dist}\big(\cC,(a/c,b/c)\big)\ll\ve/c^2$, and  hence Theorem~\ref{t9} is equivalent to \eqref{e41} for the specific type of curves $\cC$ in question.
%For instance, if $Q(x,y,z)=x^2+y^2-(rz)^2$ for some $r>0$, then $\cC$ is the circle of radius $r$ centred at the origin. In general, one can assume without loss of generality that
%$Q(x,y,1)=\widetilde Q(x-x_0,y-y_0)-1$ for some positive definite quadratic form $\widetilde Q$ of two variables.

Apparently, when attacking Problem~\ref{p5} one has to appeal to an unbounded $g_t$-orbit of $G_x\Z^3$ as opposed to bounded parts of this orbit appearing in \eqref{orbit}, where $g_t$ and $G_x$ are the same as  in \eqref{e41+}.
Indeed, assuming {that}  $\cC=\{(x,f(x):x\in U\}$ is bounded, it is a relatively simple task to verify that
\begin{equation}\label{e44}
  \eqref{e41}\quad\Longrightarrow\quad \{g_tG_x\Z^3:x\in U, t\ge0\}\text{ is unbounded in }X_3\,,
\end{equation}
while the converse requires a slight tightening of the condition on the right by replacing $U$ with any closed subset $U'$ of the interior of $U$, in which case we have that
\begin{equation}\label{e44+}
  \eqref{e41}\quad\Longleftarrow\quad \{g_tG_x\Z^3:x\in U', t\ge0\}\text{ is unbounded in }X_3\,.
\end{equation}
The argument in support of \eqref{e44} and \eqref{e44+} can be obtained on modifying the technique used for detecting rational points near manifolds that we discussed above and as detailed in any of \cite{MR2874641, MR2373145, MR2729002, divergence}. Of course, due to Margulis' theorem on {the Oppenheim} conjecture, \eqref{e41} and consequently the right hand side of \eqref{e44} hold for irrational  {ellipses}.

{Oppenheim's conjecture is only one example of problems on `{\em small values of homogeneous polynomials at integral points}'. Clearly, any problem of this ilk falls into the framework of `rational points near manifolds'. To give another example, which is of current interest and where the QnD estimate plays an important role, consider counting integral (irreducible) polynomials $P$ of degree $n$ and height $H(P)\le Q$ with relatively small discriminant $D(P)$. Indeed, for polynomials of degree $2$ the problem reduces to counting rational points near the parabola $y=x^2$, see \cite[\S2]{MR3505745}. In general, $D(P)$ can be written as a homogeneous polynomial $D(a_0,\dots,a_n)$ of the coefficients $a_0,\dots,a_n$ of $P=a_nx^n+\dots+a_0$; the degree of $D$ is $2n-2$. Thus, when $H(P)=\max\{|a_0|,\dots,|a_n|\}\le Q$, we have that $|D(P)|\ll Q^{2n-2}$. This gives rise to the following
\begin{problem}\label{p8}
Let $n\ge2$ be an integer and $v\in[0,n-1]$. Establish the asymptotic behaviour (as $Q\to\infty$) of the number $N_n(Q)$ of integral irreducible polynomials $P$ of degree $n$ and height $H(P)\le Q$ satisfying the condition
\begin{equation}\label{discr}
0<|D(P)|\ll Q^{2n-2-2v}\,.
\end{equation}
The problem can be equally restated for monic polynomials $P=x^{n+1}+a_nx^n+\dots+a_0$ of degree $n+1$.
\end{problem}
It was shown in \cite{MR3505745} that
\begin{equation}\label{vb444}
N_n(Q)\gg Q^{n+1-\frac{n+2}{n}v}
\end{equation}
for any $v\in[0,n-1]$. Quite remarkably, the proof of \eqref{vb444} represents yet another application of the QnD estimate. To be more precise, establishing \eqref{vb444} uses counting irreducible polynomials $P$ such that $P$ and its derivatives have prescribed values at points $x$ from a subset of $[-\tfrac12,\tfrac12]$ of measure at least $\tfrac12$, see \cite[Lemma~4]{MR2684299}. The latter is proved by using the QnD estimate applied to the system
$$
|P^{(i)}(x)|<\theta_i\qquad(0\le i\le n)
$$
for a suitable choice of positive parameters $\theta_i$ such that the prpoduct $\theta_0\cdots\theta_n$ is a sufficiently small constant, see \cite[Lemma~1]{MR2684299} or more generally \cite[Theorem~5.8]{MR2874641}.
In all likelihood \eqref{vb444} is sharp, but the complementary upper bound remains unknown except for $n=2$ \cite{MR3505745} and $n=3$ when $0<v<3/5$ \cite{GKK14}. Very recently, in \cite[Theorem~1.1]{DOS}, an upper bound for the number of monic irreducible polynomials of a fixed discriminant and height $H(P)\le Q$ has been established for arbitrary degrees $\ge3$. However, this recent upper bound seem to have enough room for further improvement even for monic polynomials of degree $3$ and thus finding upper bounds within Problem~\ref{p8} remains an almost entirely open challenge.
}

\section{Badly approximable points on manifolds}\label{BAP}

The notion of badly approximable points in $\R^n$ comes about by reversing the inequalities in Dirichlet's theorem with a suitably small constant. Recall again, by Dirichlet's theorem, for every $\vv y=(y_1,\dots,y_n)\in\R^n$ there are infinitely many $q\in\N$ such that
$$
\max_{1\le i\le n}|\langle qy_i\rangle|^n<q^{-1}\,,
$$
where $|\langle qy_i\rangle| $ is the distance from $qy_i$ to the nearest integer $p_i$.
Thus, the point $\vv y\in\R^n$ is \emph{badly approximable}\/ if there exists a constant $c=c(\vv y)>0$ such that
\begin{equation}\label{vb200}
\max_{1\le i\le n}\left|\langle qy_i\right\rangle| ^n\ge cq^{-1}
\end{equation}
for all $q\in\N$. More generally, given an $n$-tuple of weights
$\vv r=(r_1,\dots,r_n)\in\R^n_{\ge0}$ normalised by \eqref{weights},
the point $\vv y\in\R^n$ is called \emph{$\rr$-badly approximable} if there exists $c=c(\vv y)>0$ such that
\begin{equation}\label{e48}
\max_{1\le i\le n}|\langle qy_i\rangle| ^{1/r_i}\ge cq^{-1}
\end{equation}
for all $q\in\N$. Here, by definition, $|\langle qy_i\rangle|^{1/0}=0$. In what follows
the set of $\rr$-badly approximable points in $\R^n$ will be denoted by $\Badr$. It is a well known fact that $\Badr$ is always of Lebesgue measure zero. Therefore, in Diophantine approximation one is interested in understanding how small the sets $\Badr$ are really by using, for example, Hausdorff dimension. More sophisticated problems arise when one considers the intersections of $\Badr$ and restrictions to submanifolds of $\R^n$. This broad theme has been around for several decades and investigated in great depth, see \cite{MR0166154, MR1911218, MR2191212, MR2231044, MR2581371, MR2520102, MR0195595,
MR3425389, MR3231023, MR2846492, MR3224831, MR3284116, MR3733884, MR3457674, MR3081541} amongst many dozens of other papers on the topic. There is also a natural link, known as Dani's correspondence \cite{MR794799}, between badly approximable points in $\R^n$ and bounded orbits of the lattices
$$
\Lambda_{\vv y}=\left(\begin{array}{cc}
            I_n  & \vv y \\[1ex]
             0 & 1
           \end{array}
\right)\Z^{n+1},
$$
where $\vv y\in\R^n$ is treated as a column and $I_n$ is the identity matrix. According to Dani's correspondence, a point $\vv y\in\R^n$ is badly approximable if and only if
$g_t\Lambda_{\vv y}$ $(t\ge0)$ is bounded in the space of lattices $X_{n+1}$, where
$g_t:={\rm diag}\{e^{t},\dots,e^{t},e^{-nt}\}$. Later it was shown in \cite{MR1646538} that Dani's correspondence extends to Diophantine approximation with weights and for matrices. In particular,
a point $\vv y\in\R^n$ is $\rr$-badly approximable if and only if
{the trajectory} $g_t\Lambda_{\vv y}$ $(t\ge0)$ is bounded in $X_{n+1}$, where
$g_t:={\rm diag}\{e^{tr_1},\dots,e^{tr_n},e^{-t}\}$.

The purpose of this section is to expose the role of the QnD estimate in a recent proof given in \cite{MR3425389} that countable intersections of the sets $\Badr$ restricted to a non-degenerate submanifold of $\R^n$ have full Hausdorff dimension. The following key result of \cite{MR3425389} will be the main subject of discussion of this section.

\begin{theorem}[Theorem~1 in \cite{MR3425389}]\label{Bad:t}
Let $n,d\in\N$, $W$ be a finite or countable collection of $n$-tuples $(r_1,\dots,r_n)\in\R_{\ge0}^n$ with $r_1+\dots+r_n=1$. Assume that
\begin{equation}\label{tau}
\inf\{\tau(\rr):\rr\in W\}>0
\end{equation}
where
$$
\tau(r_1,\dots,r_n)=\min\{r_i>0:1\le i\le n\}\,.
$$
Let $\cF_n(B)$ be a finite collection of analytic non-degenerate maps defined on a ball $B\subset\R^d$. Then
\begin{equation}\label{vb1++}
\dim\bigcap_{\vv f\in\cF_n(B)}\ \bigcap_{\rr\in W}\vv f^{-1}\big(\Bad(\rr)\big)=d\,.
\end{equation}
\end{theorem}

\medskip

\noindent\textbf{Sketch of the proof} (for full details see \cite{MR3425389}).
To begin with, one uses a transference principle to reformulate $\Badr$ in terms of approximations by one linear form: $\vv y\in\R^n$ is in $\Badr$ if and only if there exists $c>0$ such that for any $H\ge1$ the only integer solution $(p,q_1,\dots,q_n)$ to the system
\begin{equation}\label{eee10}
\begin{array}{l}
|p+q_1y_1+\dots+q_ny_n|< c H^{-1},\qquad\\[1ex]
|q_i|< H^{r_i}\qquad (1\le i\le n)
\end{array}
\end{equation}
is zero, that is $p={q_1 =}\dots=q_n=0$. Another simplification is that one can assume that $d=1$, that is it suffices to deal with curves. This is due to the existence of appropriate techniques for fibering analytic non-degenerate manifolds into non-degenerate curves and Marstrand's slicing lemma, see \cite{MR3425389}.

For simplicity we will assume that $\#W=1$, $\#\cF_n(B)=1$ and $B=[0,1]$. Then \eqref{vb1++} becomes
\begin{equation}\label{S}
\dim\underbrace{\{x\in[0,1]:\vv f(x)\in\Bad(\rr)\}}_{\textstyle S}=1\,.
\end{equation}
The basic idea is to construct a Cantor set
$$
\cK:=\bigcap_{t=1}^\infty\cK_{t+m},
$$
starting from $\cK_{m}=[0,1]$, where $m$ is a large integer, and fulfilling the condition
\begin{equation}\label{e53}
\qquad\cK_{t+m}\subset\cK_{t-1+m}\smallsetminus\{x\in[0,1]:\delta(g_tu_{\vv f(x)}\Z^{n+1})< \kappa\}\qquad\text{for {$t\in\N$}}\,,
\end{equation}
where $\eta>0$ is a suitably large constant, $g_t={\rm diag}\{e^{\eta t},e^{-\eta tr_1},\dots,e^{-\eta tr_n}\}$ and $u_{\vv f(x)}$ is the same as in \eqref{new1}. By Dani's correspondence, or rather {by} its version from \cite{MR1646538}, $\cK$ is a subset of $S$ defined in \eqref{S}. The goal is thus to demonstrate that for any $\delta>0$ there exists a suitably small $\kappa>0$ such that
\begin{equation}\label{e54}
\dim \cK\ge 1-\delta.
\end{equation}

The level sets $\cK_t$ of $\cK$ are made of small `{\em building blocks}' -- closed subintervals of length $R^{-t}$ with disjoint interiors, where the parameter $R$ is a large positive integer. This requirement makes it easier to estimate the Hausdorff dimension of $\cK$. Essentially, $\cK_t$ is obtained from $\cK_{t-1}$ by chopping up each `building block' of $\cK_{t-1}$ into $R$ equal pieces and then removing some of them.
The `building blocks' that have to be removed are identified by requirement \eqref{e53}. Effectively to achieve the dimension bound in \eqref{e54} one has to show that we remove relatively little. How little is determined by a technical statement on Cantor sets originally obtained in \cite{MR2846492} and \cite{MR2838058} and developed further in \cite{MR3425389} into a notion of {\em Cantor rich} sets. Cantor rich sets are closed under countable intersections, albeit there is a mild technical condition attached to intersections, see also \cite{CWsets} for a comparison of Cantor rich sets with other similar notions. It is the nature of Cantor rich sets that allowed us to assume that $\#W=1$ and $\#\cF_n(B)=1$.

To accomplish the final goal one has to analyse the composition of the set
$$
\{x\in[0,1]:\delta(g_tu_{\vv f(x)}\Z^{n+1})< \kappa\},
$$
that is, {the set removed} in \eqref{e53}. This set is defined as the union
over all the integer points $(p,q_1,\dots,q_n)$ subject to $|q_i|< e^{\eta tr_i}$ $(1\le i\le n)$
of all the solutions $x$ to
\begin{equation}\label{km3p+}
|p+q_1f_1(x)+\dots+q_nf_n(x)|< \kappa e^{-\eta t}\,.
\end{equation}
For a fixed $(p,q_1,\dots,q_n)$ inequality \eqref{km3p+} defines a finite collection of intervals. The number of these interval is bonded by a constant depending on $n$ and $\vv f$, however the length of these intervals depends on the slope of the graph of the function $x\mapsto  p+q_1f_1(x)+\dots+q_nf_n(x)$, that is on
\begin{equation}\label{e57}
|q_1f'_1(x)+\dots+q_nf'_n(x)|
\end{equation}
and thus can vary hugely. It is convenient to combine together the intervals of similar size by sandwiching \eqref{e57} between consecutive powers of a real number. Effectively, {for some $\ell\in\Z$} one considers the system
\begin{equation}\label{km3p2}
\begin{array}{l}
|p+q_1f_1(x)+\dots+q_nf_n(x)|< \kappa e^{-\eta t},\qquad\\[1ex]
e^{\eta (\gamma t-\gamma'\ell)}\le |q_1f'_1(x)+\dots+q_nf'_n(x)|< e^{\eta\left(\gamma t-\gamma'(\ell-1)\right)},\qquad\\[1ex]
|q_i|< e^{\eta tr_i}\qquad (1\le i\le n).
\end{array}
\end{equation}
Since the maximum of \eqref{e57} is $\ll e^{\eta t\gamma}$, where $\gamma=\max\{r_1,\dots,r_n\}$, it suffices to assume that $\ell$ is non-negative. The parameter $\gamma'$ is used for convenience to eventually synchronize the (approximate) length of the intervals arising from \eqref{km3p2} with that of `building blocks' of an appropriate level of $\cK$. Indeed, for relatively small $\ell$ the intervals of $x$ arising from \eqref{km3p2} for a fixed $(p,q_1,\dots,q_n)$ are of length
\begin{equation}\label{e58}
\asymp \kappa e^{-\eta t}e^{\eta(\gamma t-\gamma'\ell)}\,.
\end{equation}
The proof uses a counting argument from the geometry of numbers to estimate the number of different points $(p,q_1,\dots,q_n)$ that give rise to a non-empty set of $x$ satisfying \eqref{km3p2} and this estimate put together with \eqref{e58} appears to be sufficient to make the Cantor rich sets work.

The problem remains in the case of relatively large $\ell$. And this is precisely the case where the QnD estimate comes to the rescue. The idea is to consider the system
\begin{equation}\label{km3p}
\begin{array}{l}
|p+q_1f_1(x)+\dots+q_nf_n(x)|< \kappa e^{-\eta t},\qquad\\[1ex]
|q_1f'_1(x)+\dots+q_nf'_n(x)|< e^{\eta t(\gamma-\ve)},\qquad\\[1ex]
|q_i|< e^{\eta tr_i}\qquad (1\le i\le n)
\end{array}
\end{equation}
where $\ve$ is a fixed constant. In practice, $\ve$ can be chosen within the limits  $\tfrac1n\le \ve\le \frac2n$. The %$x$'s satisfying
{solutions of} \eqref{km3p2} with $\ell\gg \ve t$ will fall into the set $S_t$ of solutions $x$ to \eqref{km3p}. Using the version of the QnD estimate from \cite{MR1829381} one verify that the measure of $S_t$ is
$$
\ll e^{-t\alpha\ve}\,,
$$
where $\alpha$ depends on $n$ only. In fact, if we swell the set $S_t$ up by placing a ball of radius
$$
\Delta:=\kappa e^{-\eta t}/e^{\eta t(\gamma-\ve)}
$$
around each point of $S_t$, the QnD estimate applies to this bigger set $\widehat S_t$. Due to its construction the set $\widehat S_t$ can be written as a disjoint union of intervals of length $\asymp\Delta$, while the total measure of these intervals is still $\ll e^{-t\alpha\ve}$. Hence, one gets a bound on the number of the intervals, and this bound appears good enough to complete the proof.

\medskip

{\begin{remark}
The basic idea for treating \eqref{km3p} that we described above evolved from the paper \cite{MR2069553} which deals with a very special case of Problem~\ref{p2} discussed in \S\ref{HDFR}. Indeed, the method of \cite{MR2069553}, which relies on the QnD estimate, can be easily modified to obtain the following upper bound for the Hausdorff dimension:
\begin{equation}\label{dim0}
\dim\Big\{x:~\begin{array}{l}
\text{\eqref{km3p} has a non-zero solution $(p,q_1,\dots,q_n)$}\\
\text{ for infinitely many $t\in\N$}
            \end{array}
\Big\}\le 1-c(\alpha,\ve,n)
\end{equation}
for some explicitly computable parameter $c(\alpha,\ve,n)>0$ depending only on $\alpha$, $\ve$ and $n$. Recall that within the proof of Theorem~\ref{Bad:t} the `target' set $\cK$ given by \eqref{e54} is sought to satisfy \eqref{e54} for arbitrarily small $\delta>0$. In the case of \eqref{km3p2} this goal is attained by taking $\eta$ sufficiently large and $\kappa$ suffiiciently small. Now note that the estimate \eqref{dim0} is independent of $\eta$ and $\kappa$. This means that when constructing the levels $\cK_t$ of our Cantor set, the case \eqref{km3p} `removes' a set of dimension strictly smaller than $1-\delta$ as long as we impose the condition $0<\delta<c(\alpha,\ve,n)$ with $\delta$ the same as in \eqref{e54}.
\end{remark}
}

\subsection{Further remarks}

The technical condition \eqref{tau} on the weights of approximation arises within the part of the proof of Theorem~\ref{Bad:t} that does not use the QnD. Introducing new ideas to this part, Lei Yang \cite{Yang} managed to remove \eqref{tau} completely.

Theorem~\ref{Bad:t} has a straightforward consequence {to} real numbers badly approximable by algebraic numbers. These can be defined via small values of polynomials:
$$
\cB_n=\left\{\xi\in\R:\begin{array}{l}
\exists\ c_1=c_1(\xi,n)>0\text{ such that }|P(\xi)|\ge c_1H(P)^{-n}\\
\qquad\text{for all non-zero }P\in\Z[x],\ \deg P\le n
               \end{array}
\right\}\,.
$$
As a consequence of Theorem~\ref{Bad:t} we have that for any natural number $N$ and any interval $I$ in $\R$
\begin{equation}\label{cor2}
\dim\bigcap_{n=1}^N\cB_n\cap I~=~1\,.
\end{equation}
However, Theorem~\ref{Bad:t} leaves the following problem open: {\em show that \eqref{cor2} holds when $N=\infty$.}
The generalisation of Yang \cite{Yang} that removes condition \eqref{tau} does not solve this problem. However, most recently, it has been resolved in \cite{BNY1} on showing that the sets
$$
\{x\in\R:(x,\dots,x^n)\text{ is badly approximable}\}
$$
are winning. Previously, this was shown in dimension $n=2$ \cite{MR3733884}. More generally, it is shown in \cite{BNY1} that for any $n\in\N$ and any $n$-tuple $\vv r$ of weights the set of $\vv r$-badly approximable points on any non-degenerate analytic curve in $\R^n$ is absolute winning. We note that the results of \cite{BNY1} represent yet another powerful application of the QnD, this time for fractal measures as established in \cite{MR2134453}. Another remarkable application of the QnD for fractal measures is the proof that the sets $\Bad(\vv r)$ are hyperplane absolute winning established in \cite{BNY2}.

We note that the above exposition is not a complete account of known applications of the QnD estimates. There is no doubt that many new exciting applications are still awaiting to be discovered!

{\small\def\cprime{$'$} \def\polhk#1{\setbox0=\hbox{#1}{\ooalign{\hidewidth
  \lower1.5ex\hbox{`}\hidewidth\crcr\unhbox0}}} \def\cprime{$'$}
  \def\cprime{$'$}

}

%{\small
%%---- plain!, alpha, unsrt, abbrv, siam!, amsalpha, apalike, acm
%\bibliographystyle{alpha}
%\bibliography{%
%C:/DATA/Vitia/Papers/BIBLIOGRAPHY/arg2018,%
%C:/DATA/Vitia/Papers/BIBLIOGRAPHY/beresnevich,%
%C:/DATA/Vitia/Papers/BIBLIOGRAPHY/bernik,%
%C:/DATA/Vitia/Papers/BIBLIOGRAPHY/kleinbock,%
%C:/DATA/Vitia/Papers/BIBLIOGRAPHY/dani,%
%C:/DATA/Vitia/Papers/BIBLIOGRAPHY/ghosh,%
%C:/DATA/Vitia/Papers/BIBLIOGRAPHY/margulis,%
%C:/DATA/Vitia/Papers/BIBLIOGRAPHY/collected,%
%C:/DATA/Vitia/Papers/BIBLIOGRAPHY/Diophant}
%}

\medskip
\medskip
\medskip

{\footnotesize

\begin{minipage}{0.9\textwidth}
\footnotesize V. Beresnevich\\
Department of Mathematics, University of York, Heslington, York, YO10 5DD, England\\
{\it E-mail address}\,:~~ \verb|victor.beresnevich@york.ac.uk|\\
\end{minipage}

\begin{minipage}{0.9\textwidth}
\footnotesize D. Kleinbock\\
{\it E-mail address}\,:~~ \verb|kleinboc@brandeis.edu|\\
\end{minipage}

}

\end{document}